\documentclass[12pt]{amsart}

\usepackage{amssymb,amscd,epsfig,subfigure,times}

\numberwithin{equation}{section}
\newtheorem{Theorem} {Theorem}
\newtheorem* {Proposition} {Proposition}
\newtheorem {Lemma} {Lemma}
\theoremstyle{remark}
\newtheorem{Remark} {Remark}

\def \C {{\mathbb C}}
\def \R {{\mathbb R}}
\def \Q {{\mathbb Q}}
\def \Z {{\mathbb Z}}

\begin{document}

\title{Vortex dynamics on a cylinder}

\author{James Montaldi}\author{Anik Souli\`ere}\author{Tadashi Tokieda}

\address{Department of Mathematics, UMIST, PO Box 88, Manchester 
M60 1QD, UK}
\email{j.montaldi@umist.ac.uk}

\address{D\'ept.~de Math\'ematiques, Universit\'e de
Montr\'eal, C.P.~6128, succ.~Centre-Ville, Montr\'eal H3C 3J7,
Canada}
\email{souliere@dms.umontreal.ca}
\email{tokieda@dms.umontreal.ca}

\date{\today}

\begin{abstract}
  Point vortices on a cylinder (periodic strip) are studied
  geometrically.  The Hamiltonian formalism is developed, a
  non-existence theorem for relative equilibria is proved, equilibria
  are classified when all vorticities have the same sign, and several
  results on relative periodic orbits are established, including as
  corollaries classical results on vortex streets and leapfrogging.
\end{abstract}

\maketitle

\section{Introduction}

Spatially periodic rows of point vortices in a 2-dimensional ideal
fluid have long attracted the attention of fluid dynamicists, one of
the earliest and the most popular instances being K\'arm\'an's vortex
street \cite{Karman}, \cite[photos~94--98]{VanDyke}.  The general
problem is as follows: analyse the motion of an infinite configuration
consisting of vortices $z_1, \ldots, z_N\in \C$ with vorticities
$\Gamma_1, \ldots, \Gamma_N\in \R$ together with their translates $\{
z_k + 2\pi rm \> | \> k = 1, \ldots, N, m\in\Z \}$, where $2\pi r > 0$
is the spatial period of translation.  Traditionally the problem is
analysed on the plane $\C$, but in this paper we place the vortices on
a cylinder $\C/2\pi r\Z$ (fig.~1).
\begin{figure}[ht]
  \begin{center}
    \includegraphics[width=9cm]{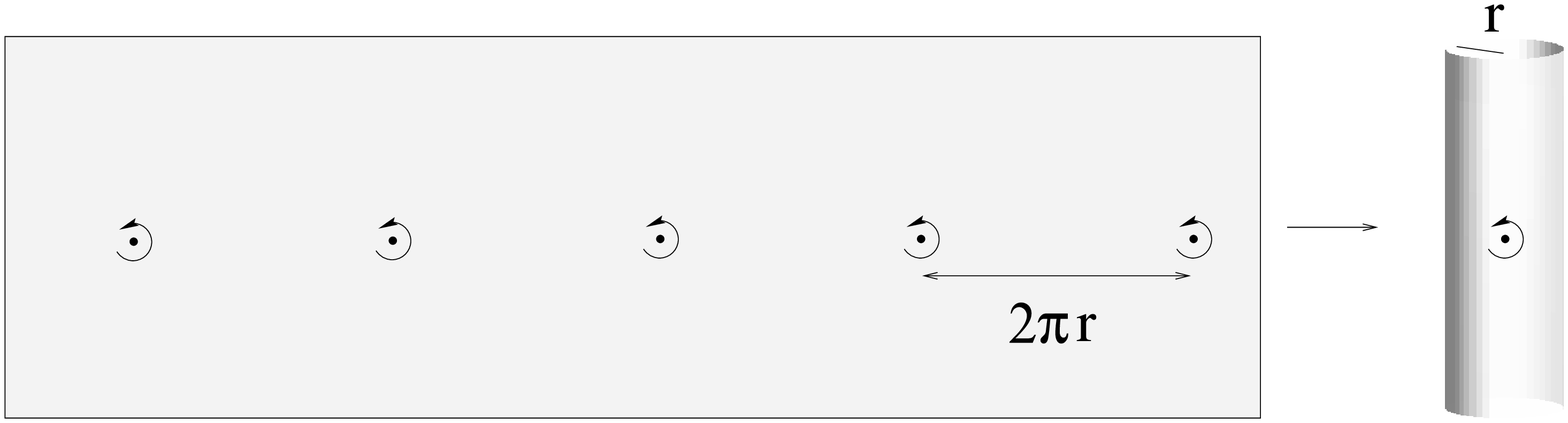}
    \caption{}
  \end{center}
\end{figure}
\noindent Though the two pictures---periodic planar and 
cylindrical---are for most purposes equivalent, as we shall see there
are advantages, both conceptual and computational, to working on a
cylinder rather than on the plane.  The proviso `for most purposes' is
necessary because the cylindrical picture posits that everything in
the dynamics be $2\pi r$-periodic, whereas in the planar picture one
could allow, for example, non-periodic perturbations to the periodic
row.  Physically, however, perturbations are usually due to some small
change in the mechanism generating the vortex row, and such a change
generates spatially periodic perturbations.  Therefore it is natural
to look at the cylindrical picture first.

We shall be interested in how vortices move relative to one another,
more precisely in their dynamics modulo the translational action of
the symmetry group $\C/2\pi\Z$.  The basic objects of interest are
relative equilibria and relative periodic orbits.  A {\it relative
  equilibrium\/} is a motion of vortices that lies entirely in a group
orbit (i.e.\ it looks stationary up to translation), and a {\it
  relative periodic orbit\/} is a motion that revisits the same group
orbit after some time (i.e.\ it looks periodic in time up to
translation).  Equilibria and periodic orbits in the ordinary sense
are special examples of relative equilibria and relative periodic
orbits.  When we wish to exclude ordinary equilibria or periodic
orbits, we speak of relative equilibria or relative periodic orbits
{\it with nonzero drift}.

As on the plane, dynamics of point vortices on a cylinder lends itself
to a Hamiltonian formalism.  The model presented here is then a
finite-dimensional Hamiltonian approximation to the vortex dynamics of
the Euler equation.  This approximation is mathematically very rich and
in the context of the plane can claim a pedigreed history
\cite[chap.~VII]{Lamb}, \cite{Villat}, \cite{Aref}.  Conversely, 
the motion of point
vortices is amenable to desingularization to a solution of the Euler
equation.

For vortices on the plane or on a sphere, an extensive theory of
relative equilibria is available (especially when the vorticities are
identical) \cite{Arefetal}, \cite{LMR}.  In contrast, apart from a
study on 3 vortices \cite{ArefStremler}, no literature seems to exist
on relative equilibria and relative periodic orbits of $N$ vortices on
a cylinder.  In this paper we develop the Hamiltonian formalism for
vortex dynamics on a cylinder (section~2), prove that if the
vorticities do not sum to zero a cylinder supports no relative
equilibrium with nonzero drift (section~3), classify equilibria when
all vorticities have the same sign (section~3), show that 3 vortices
form a relative periodic orbit for `small' initial conditions or for
vorticities dependent over $\Q$ with zero sum, and establish several
results on a class of relative periodic orbits called leapfrogging
\cite[photo~79]{VanDyke} (section~4), which may be regarded as
splitting of K\'arm\'an's vortex street.

   If the vorticities have nonzero sum, the action of the symmetry 
group $\C/2\pi\Z$ does not have 2 globally defined first integrals 
(conserved quantities) associated to it: the subgroup of horizontal
translations $\R/2\pi\Z$ has a first integral, but not the 
subgroup of vertical translations $i\R$.  One of the novelties of 
the present work is to exploit {\it local\/} first integrals 
(Theorems 2, 3, 4).

Many of the results have analogues in the theory of vortices on a
torus, i.e.\ for spatially biperiodic arrays of vortices.

\section{Hamiltonian formalism of vortices on a cylinder}

Throughout the paper {\it cylinder\/} means the surface $\C/2\pi r\Z
\simeq (\R/2\pi r)\times \R$, where $r > 0$ is some fixed constant,
the {\it radius\/} of the cylinder.  The coordinate $z = x + iy$ on
$\C/2\pi r\Z$ is to be read modulo $2\pi r$, i.e.\ $x\equiv x + 2\pi
rn$ for all $n\in \Z$~; the $x$-axis (which is a circle) is {\it
  horizontal}, the $y$-axis {\it vertical}.  The phase space for the
motion of vortices $z_1, \ldots, z_N$ with vorticities $\Gamma_1,
\ldots, \Gamma_N$ is the product of $N$ copies of the cylinder with
diagonals removed (to exclude collisions).  The Hamiltonian is a
weighted combination $H(z_1, \ldots, z_N) = \sum_{k < l}
\Gamma_k\Gamma_l \psi(z_k, z_l)$ of Green's function $\psi$ for the
Laplacian on the cylinder: $\nabla^2\psi(z, z_0) = -\delta_{z_0}(z)$
(see e.g.\ \cite[section~2]{Tokieda}).  Hamilton's equations are
$$
\frac{dz_k}{dt} = \frac{2}{i}\frac{\partial
  H}{\partial(\Gamma_k\overline{z}_k)}, \qquad (k = 1, \ldots, N).
$$
The quickest way to derive the Hamiltonian on a cylinder is to
periodize Green's function on the plane $\psi(z_k, z_l) = -
\frac{1}{2\pi}\log|z_k - z_l|$ by taking into account contributions
from $2\pi\Z$-translates.  Formally the periodized Hamiltonian becomes
$$
-\frac{1}{2\pi}\sum_{n \in \Z}\sum_{k < l} \Gamma_k\Gamma_l\log|z_k
- z_l - 2\pi rn|,
$$
which, as it stands, diverges.  But since additive constants in $H$
do not affect the dynamics, we can subtract off a constant divergent
series to force the remaining functional part to converge.  
Jettisoning $-\frac{1}{2\pi}\sum_n
\sum_{k<l}\Gamma_k\Gamma_l\log|2\pi rn|$ and pairing terms in $n$ and
$-n$,
\begin{equation}\label{complexh}
H = -\frac{1}{2\pi}
    \sum_{k < l} \Gamma_k\Gamma_l
    \log \left| 
    (z_k - z_l)
    \prod_{n \geqslant 1}
    \left( 1-\left(\frac{z_k - z_l}{2\pi rn}\right)^{\!\! 2}\right)
    \right|
  = -\frac{1}{2\pi}
    \sum_{k < l} \Gamma_k\Gamma_l
    \log\left|\sin\frac{z_k - z_l}{2r}\right|.
\end{equation}
The equations of motion on a cylinder are therefore
\begin{equation}\label{complexeq}
\frac{dz_k}{dt} = \frac{i}{4\pi r}
   \sum_{l, l\neq k} 
   \Gamma_l\, {\rm cotan\,}\frac{\overline{z_k} - \overline{z_l}}{2r},
\qquad (k = 1, \ldots, N).
\end{equation}
\noindent
For reference, we list expressions in real coordinates:
\begin{equation}\label{realh}
H = -\frac{1}{4\pi}\sum_{k<l} \Gamma_k\Gamma_l
    \log\left\{ \sin^2 \left( \frac{x_k - x_l}{2r}\right) +
             {\rm sinh}^2\left(\frac{y_k - y_l}{2r}\right) \right\},
\end{equation}
\begin{equation}\label{realeq}
\begin{cases}
  \dfrac{dx_k}{dt} = -&\displaystyle \dfrac{1}{8\pi r}\sum_{l, l\neq
    k}\Gamma_l \dfrac{{\rm sinh\,}\dfrac{y_k - y_l}{r}}
  {\sin^2\left(\dfrac{x_k - x_l}{2r} \right) +
    {\rm sinh}^2\left(\dfrac{y_k - y_l}{2r} \right)} \\[36pt]
  \dfrac{dy_k}{dt} = &\displaystyle \dfrac{1}{8\pi r}\sum_{l, l\neq
    k}\Gamma_l \dfrac{\sin\dfrac{x_k - x_l}{r}}
  {\sin^2\left(\dfrac{x_k - x_l}{2r} \right) + {\rm
      sinh}^2\left(\dfrac{y_k - y_l}{2r} \right)}
\end{cases},\qquad (k = 1, \ldots, N).
\end{equation}

One noteworthy feature of (\ref{realeq}) is that as $y_k - y_l\to
\infty$ (infinite vertical separation), the velocity induced by $z_l$
on the vortex $z_k$ does not decay to $0$, but tends to $\Gamma_l/4\pi
r$, as is obvious upon calculating in the planar theory the
circulation around a tall window of width $2\pi r$ enclosing $z_l$.
Another way to interpret the feature is to note that in the planar
theory, up to rescaling, stretching vertical separation amounts to
narrowing the spatial period $2\pi r\to 0$~; the latter limit produces
a {\it vortex sheet\/} (or more aptly {\it vortex line\/} in this
2-dimensional theory), which induces a velocity field constant above
(and the opposite constant below) the sheet independently of the
distance to the sheet.  This is exactly as in 2-dimensional
electromagnetism or gravity where the force induced by a homogeneous
charge or mass distribution along an infinite line is independent of
the distance to the line.

   Physically, periodizing the plane with period $2\pi r$ and
considering $N$ vortices on the resulting cylinder is the same as
periodizing with period $2\pi rn$ and considering $nN$ vortices on the
resulting wider cylinder.  The equivalence between these
periodizations is trivial yet sometimes useful:

\begin{Proposition}
   Let $z_1, \ldots, z_N$ be vortices with vorticities $\Gamma_1, 
\ldots, \Gamma_N$ on a cylinder of radius $r$.
Next let $z_1, \ldots, z_N, z_1+2\pi r, \ldots, z_N+2\pi r, \ldots,
z_1+2\pi rn, \ldots, z_N+2\pi rn$ be their `$n$-fold copies' with
corresponding vorticities on a cylinder of radius $rn$, where $n$ is
any strictly positive integer.  Then the dynamics on the cylinder of
radius $rn$ covers the dynamics on the cylinder of radius $r$.
\end{Proposition}

In particular, given a relative equilibrium or a relative periodic
orbit, we can reel
off infinite families of relative equilibria or relative periodic
orbits at no extra cost by replicating the configuration sideways on a wider
cylinder.

\begin{Remark}
  A torus has the form $\C/(\pi\Z + \tau\pi\Z)$, where the parameter
  $\tau\in \C$, ${\rm Im}\tau > 0$ controls the conformal class.  The
  Hamiltonian is
  $$
  H = -\frac{1}{2\pi} \sum_{k < l} \Gamma_k\Gamma_l \left\{
    \log\left| \vartheta_1(z_k - z_l|\tau)\right| - \frac{({\rm
        Im}(z_k - z_l))^2}{\pi {\rm Im}\tau } \right\},
  $$
  where $\vartheta_1$ is the 1st Jacobian theta function
  \cite{O'Neil}, \cite{StremlerAref}, \cite{Tokieda}.
\end{Remark}

A cylinder has a translational symmetry of $\C/2\pi r\Z$ acting on
itself, hence acting diagonally on the phase space.  The plane has a
supplementary rotational symmetry $z\mapsto e^{i\theta}z, \theta\in
\R$~; this is lost on the cylinder.  Via Noether's theorem the
translational symmetry of $\C/2\pi r\Z$ should give rise to a first
integral, a {\it momentum map\/} $(z_1, \ldots, z_N) \mapsto \sum_k
\Gamma_k z_k$, but there is a rub: because $z$'s are defined only
modulo $2\pi r$ this `momentum map' is not well-defined as a map to
the dual of the Lie algebra of the symmetry group $\C/2\pi\Z$.  Nor is
it advisable to treat this `momentum map' as a multi-valued function,
for generically $\Gamma_1, \ldots, \Gamma_N$ are independent over $\Q$
and so the ambiguity $\{ 2\pi r\sum_k \Gamma_k n_k\> |\> n_1, \ldots,
n_N\in \Z \}$ in the value of the `map' is dense in $\R$.
Nevertheless, the momentum map is {\it locally\/} (i.e.~on each chart)
well-defined.  From now on, whenever we write $\sum_k \Gamma_k z_k$,
some suitable ad hoc chart will be understood.

When $\sum_k \Gamma_k \neq 0$, the {\it center of vorticity\/} $\sum_k
\Gamma_k z_k/\sum_k \Gamma_k$ is a more intuitive first integral
\cite[art.~154]{Lamb}.  The next result provides a substitute for
center of vorticity when $\sum_k \Gamma_k = 0$.
\begin{Theorem}
  Let $\{ z \}$ be vortices on the plane or on a cylinder whose
  vorticities sum to zero: $\sum \Gamma = 0$.  Suppose the vortices
  are partitioned into two groups $\{ z'\}$, $\{ z''\}$ and within
  each group $\sum \Gamma' \neq 0$, $\sum \Gamma'' \neq 0$, so that
  the center of vorticity for each group is well-defined.  Then the
  vector connecting the two centers of vorticity is a local first
  integral (fig.~2).
\end{Theorem}
\begin{figure}[ht]
  \begin{center}
    \includegraphics[width=3cm]{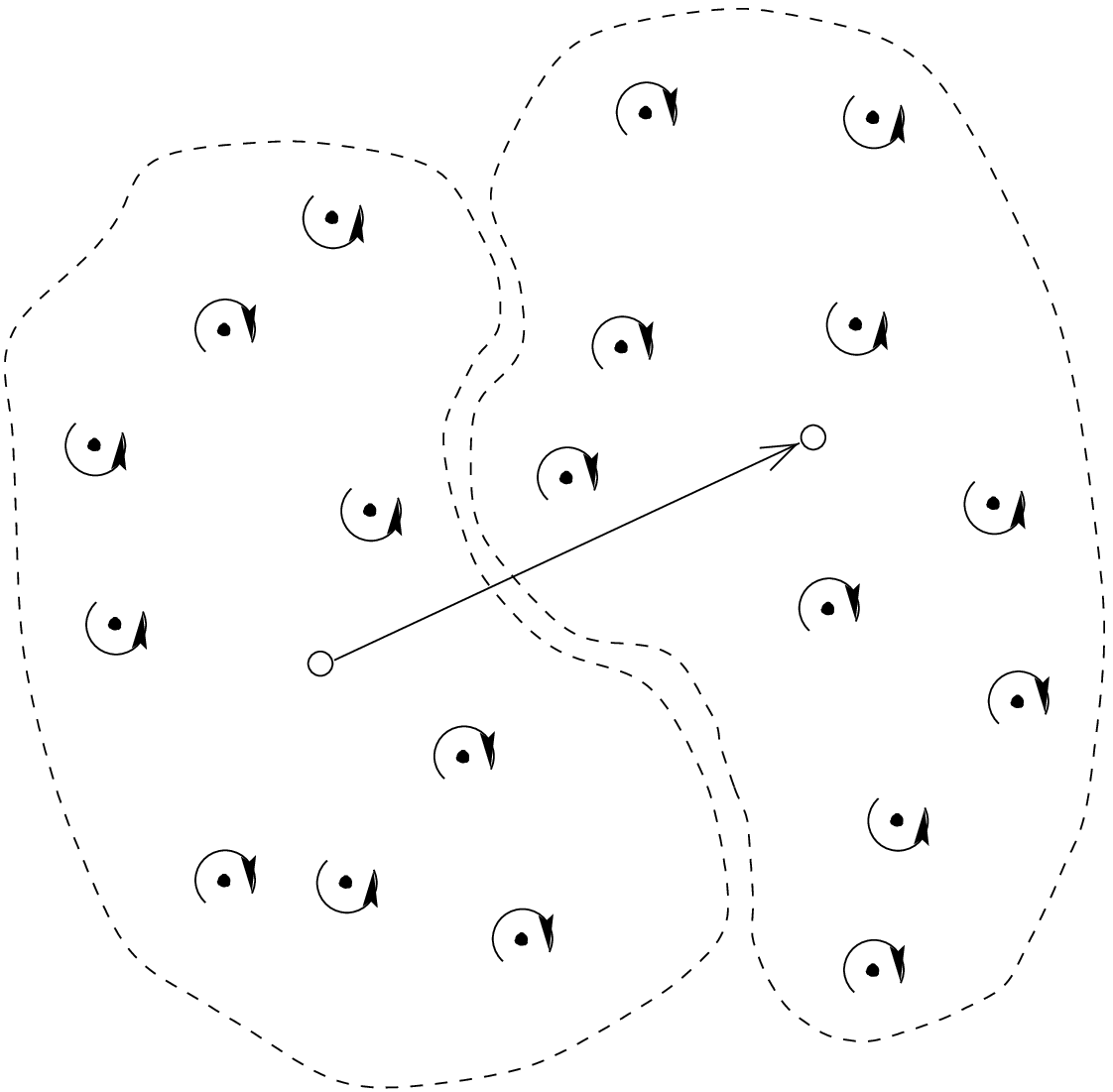}
    \caption{}
  \end{center}
\end{figure}
\begin{proof}
  Since $\sum \Gamma' + \sum \Gamma'' = 0$, the vector in question is
  $$
  \frac{\sum \Gamma' z'}{\sum \Gamma'} - \frac{\sum \Gamma''
    z''}{\sum \Gamma''} = \frac{\sum \Gamma' z'}{\sum \Gamma'} +
  \frac{\sum \Gamma'' z''}{\sum \Gamma'} = \frac{\sum \Gamma z}{\sum
    \Gamma'},
  $$
  and $\sum \Gamma z$ is a local first integral.
\end{proof}

Theorem~1 is serviceable in many problems.  The simplest illustration
is the motion of a {\it vortex pair\/} $z_1, z_2$ with vorticities
$\Gamma, -\Gamma$ \cite[photos~77, 78]{VanDyke}.  Treating $z_1$ as
one group and $z_2$ as the other group, we check against Theorem~1
that $z_2 - z_1$ is constant during the motion.  In fact, according to
(\ref{realeq}) the vortex pair on a cylinder forms a relative
equilibrium moving with slope
$$
- \sin\dfrac{x_2 - x_1}{r}\Bigr/{\rm sinh\,}\dfrac{y_2 - y_1}{r}.
$$
When $x_2 - x_1 = 0$ or $\pi r$ the pair moves horizontally: the
corresponding configurations on the plane are the unstaggered or fully
staggered cases of K\'arm\'an's vortex street.  When $z_1, z_2$ are in
general position, the corresponding vortex street on the plane
translates at an angle to the horizontal, a case studied in
\cite{Maue}.  The `plane limit' $r\to \infty$ yields the angle of
progression of a vortex pair on the plane $-(x_2 - x_1)/(y_2 - y_1)$.
For a beautiful study of the stability of variants of vortex streets,
see \cite{Imai}.

\section{Relative equilibria}

The first fact about relative equilibria of vortices on a cylinder is
that there are not many of them.

\begin{Theorem}
  Let $z_1, \ldots, z_N$ be vortices with vorticities $\Gamma_1,
  \ldots, \Gamma_N$ on a cylinder $\C/2\pi r\Z$. Suppose $\sum_k
  \Gamma_k \neq 0$.  Then all relative equilibria are in fact
  equilibria.  Moreover, if all $\Gamma$'s have the same sign, then
  for each cyclic ordering there exists a unique (up to translation by
  $\C/2\pi r\Z$) equilibrium, and all the vortices are aligned on a
  single horizontal circle.
\end{Theorem}

\begin{proof}
  If $z_1, \ldots, z_N$ form a relative equilibrium, then all $z$'s
  move with some common drift velocity $v$.  The local first integral
  should not vary:
  $$
  0 = \frac{d}{dt}\sum_k \Gamma_k z_k = v\sum_k \Gamma_k,
  $$
  so $\sum_k \Gamma_k = 0$ or else $v = 0$.
  
  If the vortices are not aligned on a single horizontal circle, pick
  a `top vortex' (one with maximal $y$-coordinate) and a `bottom
  vortex' (one with minimal $y$-coordinate).  If all $\Gamma$'s
  have the same sign, then by (\ref{realeq}) the velocities of the top
  and bottom vortices must have $x$-components with opposite signs, so
  this position cannot constitute an equilibrium.
  
  Now suppose all the vorticities are of the same sign.  Fix a cyclic
  ordering of the vortices, and place the vortices in order on a
  single horizontal circle.  The Hamiltonian is given by
  $$
  H = -\frac{1}{4\pi}\sum_{k < l} \Gamma_k\Gamma_l \log \sin^2
  \left(\frac{x_k - x_l}{2r}\right).
  $$
  This is a convex function of $x_1, \ldots, x_N$ by the same
  argument as for \cite[Theorem~4.8]{LMR}: one first checks that the
  second derivatives satisfy $\partial^2 H/\partial x_k\partial x_l <
  0$ for $k\neq l$ and $\partial^2 H/\partial x_k^2 > 0$, $\sum_l
  \partial^2 H/\partial x_k\partial x_l = 0$ for each $k$~; it then
  follows from a variant of Gershgorin's theorem (Lemma
  \ref{Gershgorin} below) that 0 is a simple eigenvalue of the Hessian
  of $H$ and all other eigenvalues are strictly positive.
  Consequently on each connected component of the domain of definition
  there is a unique minimum and no other critical point, and different
  connected components correspond to different cyclic orderings.
\end{proof}

\begin{Lemma} \label{Gershgorin}
  Let $A = (a_{kl})$ be a symmetric $N\times N$ matrix satisfying
$a_{kl} < 0$ for $k\neq l$, and $a_{kk} > 0$, 
$\sum_{l=1}^N a_{kl} = 0$ for each $k$.
Then $0$ is a simple eigenvalue of $A$ and all other eigenvalues are
  strictly positive.
\end{Lemma}

\begin{proof}
  Let $u = (u_1, \ldots, u_n)^T$ be an eigenvector of $A$ with eigenvalue 
  $\lambda$, normalized so that there is an index $k$ for which
$u_k = 1$ and $|u_l| \leqslant 1$ for all $l$.  The $k$th row of 
the equation $Au = \lambda u$ is  
$a_{kk} + \sum_{l, l\neq k} a_{kl}u_l =\lambda$, which in view of the
hypotheses on $a_{kl}$ may be written
$\sum_l |a_{kl}|(1 - u_l) =\lambda$.
But $1 - u_l \geqslant 0$ and $|a_{kl}| > 0$ for each $l$; it 
follows that $\lambda \geqslant 0$ and $\lambda = 0$ if and only 
if all $u_l = 1$.  On the other hand, $(1, \ldots, 1)^T$ is 
obviously an eigenvector with eigenvalue 0. 
\end{proof}

If the vortices are placed on a single horizontal circle so that
successive vorticities have alternating signs, then we also get the
existence of an equilibrium, though the uniqueness problem is open as
the function is no longer convex.  In full generality, if the signs
are neither the same nor alternating, the argument for existence fails
as $H\to +\infty$ for some collisions and $\to -\infty$ for others.

\begin{Remark}
  For $N = 2$, if $\Gamma_1 + \Gamma_2 \neq 0$, we have generically a
  periodic orbit and exceptionally an equilibrium of antipodal
  vortices $z, z + \pi r$ or a separatrix connecting equilibria.  For
  $N > 2$, if $\sum_k \Gamma_k \neq 0$ but $\Gamma$'s have mixed
  signs, equilibria are less severely constrained.  For example, for
  $N = 3$, let $z_1, z_2$ be vortices with vorticities $\Gamma_1,
  \Gamma_2 > 0$.  To secure an equilibrium, the third vortex $z_3$
  with vorticity $\Gamma_3 < 0$ must be placed at one of the 2
  stagnation points of the velocity field induced by $z_1, z_2$, given
  in view of (\ref{complexeq}) as roots of
$$
  \Gamma_1\, {\rm cotan\,}\frac{z - z_1}{2} + \Gamma_2\, {\rm
    cotan\,}\frac{z - z_2}{2} = 0.
  $$
  Having chosen $z_3$ as one of the roots and thereby immobilized
  $z_3$, adjust $\Gamma_3$ so as to immobilize $z_1$~:
  $$
  \Gamma_2\, {\rm cotan\,}\frac{z_1 - z_2}{2} + \Gamma_3\, {\rm
    cotan\,}\frac{z_1 - z_3}{2} = 0.
$$
  Then $z_2$ too is automatically immobilized:
  $$
  \Gamma_3\, {\rm cotan\,}\frac{z_2 - z_3}{2} + \Gamma_1\, {\rm
    cotan\,}\frac{z_2 - z_1}{2} = 0.
  $$
  The upshot is that given any $z_1, z_2$ with vorticities of the
  same sign, we have 2 positions to place $z_3$ with the right
  vorticity of the opposite sign to secure an equilibrium.  For
  example, vortices $z_1, z_2$ both of vorticity $\Gamma$ such that
  $z_2 - z_1 = 2ib$ are immobilized by the adjunction of a vortex
  $(z_1 + z_2)/2$ of vorticity
  $$
  \Gamma \left( \frac{1}{2} {\rm sech}^2\frac{b}{2r} - 1 \right).
  $$
  This is always less than $ -\Gamma/2$ and in the plane limit
  $r\to \infty$ tends to the corresponding value in the planar theory
  $-\Gamma/2$.  On the other hand, in the `vortex sheet limit' $b\to
  \infty$ this tends to $-\Gamma$, also as it should.  Similarly,
  vortices $z_1, z_2$ of vorticity $\Gamma$ such that $z_2 - z_1 = 2a$
  are immobilized by the adjunction of a vortex $(z_1 + z_2)/2$ of
  vorticity
  $$
  \Gamma \left( \frac{1}{2} {\rm sec}^2\frac{a}{2r} - 1 \right).
  $$
  In the planar limit this tends again to $-\Gamma/2$.  On the
  other hand, it is $0$ when $a = \pi r/2$~: $z_1, z_2$ are antipodal
  on the cylinder and are stationary already by themselves.  When
  $a\to \pi r$, $z_1, z_2$ nearly meet at the back and a stronger
  and stronger vortex is required at the front to prevent them from
  moving.
\end{Remark}

\begin{Remark}
  Now suppose $\sum_k \Gamma_k = 0$.  It was pointed out at the end of
  section~2 that a vortex pair $N = 2$ is always a relative
  equilibrium.  For $N = 3$, Aref and Stremler \cite{ArefStremler}
  made a detailed study of relative equilibria; the patterns of some
  trajectories are surprisingly complicated.  For $N > 3$ and $N$ even,
  we have for any $a, b > 0$ a family of relative
  equilibria consisting of $n = N/2$ vortices with vorticity $\Gamma$ at
\begin{equation}\label{uppertrain}
ib, ib + \frac{2\pi r}{n}, \ldots, ib + (n-1)\frac{2\pi r}{n},
\end{equation}
and $n$ vortices with vorticity $-\Gamma$ at
\begin{equation}\label{lowertrain}
a - ib, a - ib + \frac{2\pi r}{n}, \ldots, 
a - ib + (n-1)\frac{2\pi r}{n}.
\end{equation}
This is merely a crowded vortex street with spatial period $2\pi r/n$,
or equivalently a single vortex pair on a thinner cylinder of radius
$r/n$ (see stability calculations in \cite{Domm}).  No essentially
different family of relative equilibria seems to be known for $N > 3$.

Incidentally, even the trivial equivalence between 1 vortex on a
cylinder of radius $r$ and $n$ horizontally equidistributed vortices
on a cylinder of radius $nr$ leads to amusing identities
\cite{Arefetal}~: for example, equating the induced velocity fields and
rescaling the variables in (\ref{complexeq}),
$$
\frac{1}{n}\sum_{l = 1}^n {\rm cotan\,}\frac{z + \pi l}{n} = {\rm
  cotan\,}z, \quad \forall z\in \C.
$$
\end{Remark}

\begin{Remark}
  On the plane equilibria do not exist either when all $\Gamma$'s are
  of the same sign (even the possibility of a horizontal circle is
  lost), and the non-existence of translational relative equilibria
  with nonzero drift when $\sum_k \Gamma_k \neq 0$ holds also on the
  plane and on a torus; the proof carries over verbatim from the
  cylindrical theorem.  A torus, however, accommodates more varied
  families of equilibria: for example, $n_1n_2$ vortices with
  identical vorticity $\Gamma$ placed on a sub-lattice $(\pi/n_1)\Z +
  (\tau\pi/n_2)\Z$ form an equilibrium \cite{Tokieda}.  Many further
  patterns of equilibria may be designed on a torus with identical or
  alternating vortices.
\end{Remark}

\section{Relative periodic orbits}

Once a relative equilibrium of vortices is known, a frequently
successful recipe for creating relative periodic orbits consists in
{\it splitting\/} the vortices.  Assume the vortices $z_1, \ldots,
z_N$ with vorticities $\Gamma_1, \ldots, \Gamma_N$ form a relative
equilibrium.  Let us split each $z_k$ into a cluster, near the
original position of $z_k$, of $n_k$ vortices $z_{k,1}, \ldots,
z_{k,n_k}$ whose vorticities are of the same sign and sum to
$\Gamma_k$.  We expect the child vortices $z_{k,1}, \ldots, z_{k,n_k}$
to orbit around one another and remain a cluster, while seen from far
away they still look like the original parent vortex $z_k$ with
vorticity $\Gamma_k$.  It is reasonable to conjecture that for
suitable intial configurations the child vortices form a relative
periodic orbit, and for perhaps generic splittings they form a relative
{\it quasi-periodic\/} orbit.

A vortex pair on a cylinder, which corresponds in the planar picture
to K\'arm\'an's vortex street, is a relative equilibrium.  In this
section we shall create various relative periodic orbits by splitting
a vortex pair; as a special case we recover the phenomenon classically
known in the planar picture as leapfrogging. In Theorem~3 we split one
of the vortices, while in Theorem~4 we split both. The split is
measured by a complex variable $\zeta=\xi+i\eta$ (or rather by
$2\zeta$), and we are principally interested in small values of
$|\zeta|$.  In all the formulae the radius of the cylinder is
normalized to $r = 1$; denormalization is a matter of dimensional
analysis.  Later in the section additional classes of relative
periodic orbits are described.

Take a vortex pair at $c, -c$, where $c = a + ib\in \C$.  We split it
into 3 or 4 vortices as in fig.~3: the left diagram illustrates
Theorem~3; the middle one Theorem 4, case $-b(1 +
\Gamma/\Gamma')/2 < \eta < b(1 + \Gamma'/\Gamma)/2$~; the right one
case $b(1 + \Gamma'/\Gamma)/2 < \eta$.  Theorem~4, case $\eta <
-b(1 + \Gamma/\Gamma')/2$ is like the right diagram reflected
laterally with $\Gamma, \Gamma'$ interchanged.  

\begin{figure}[ht]
  \begin{center}
    \includegraphics[width=12cm]{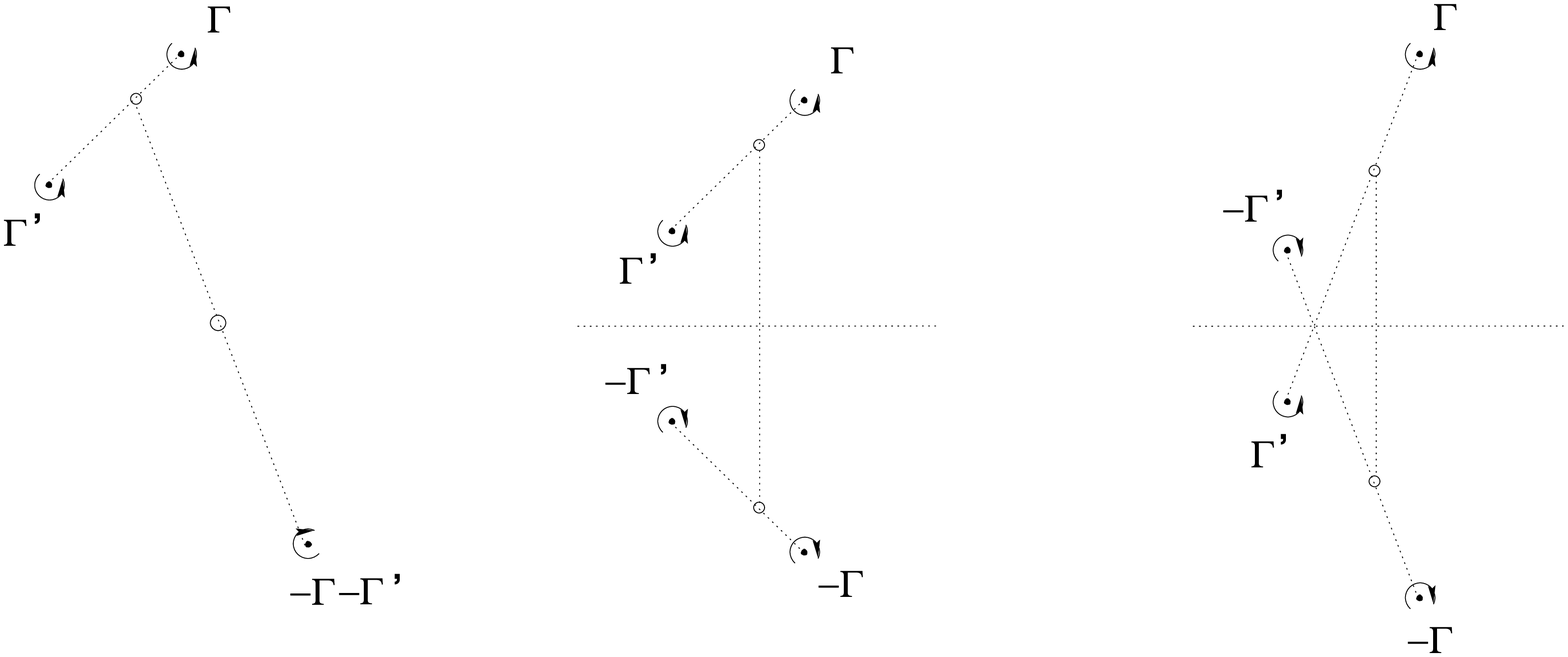}
    \caption{}
  \end{center}
\end{figure}

\begin{Theorem}
  Let $c\in \C\diagdown \{ 0 \}$.  On a cylinder, consider the
  configuration of\/ $3$ vortices with vorticities $\Gamma, \Gamma',
  -\Gamma-\Gamma'$ (\/$\Gamma$ and $\Gamma'$ being of the same sign)
  at
  $$
  c + \frac{2\Gamma'}{\Gamma + \Gamma'}\zeta,\quad c -
  \frac{2\Gamma}{\Gamma + \Gamma'}\zeta,\quad -c.
  $$
  There exists an open punctured neighborhood of $\zeta = 0$ such
that for every initial condition $\zeta(0)\neq 0$ in this
neighborhood, these vortices form a relative periodic orbit.  
If $\Gamma/\Gamma'\in
  \Q$, then for a generic choice of $\zeta(0)$ (no restriction on its
  size) these vortices form a relative periodic orbit, and for
  isolated choices of $\zeta(0)$ they form a relative equilibrium or
  a separatrix connecting relative equilibria.
\end{Theorem}

\noindent 
Combined with Proposition of section 2, Theorem 3 gives relative 
equilibria and relative periodic orbits of $N = 3n$ vortices for
all $n \geqslant 1$.
The result for $N = 3$ when $\Gamma/\Gamma' \in \Q$ 
is in \cite{ArefStremler}, but we give a somewhat different proof.
The relative periodicity for small $\zeta(0)$ is new.

The proof invokes the following elementary lemma.

\begin{Lemma} \label{lemma:punctures}
  Let $H$ be a function with only nondegenerate critical points on a
  compact surface with $p$ punctures such that $|H|\to \infty$ near
  each puncture.  Then the generic level sets of $H$ are disjoint
  unions of loops.  If $p > 2$, then besides loops there exist
  isolated saddles and separatrices connecting the saddles.
\end{Lemma}

\noindent The idea now is to use symmetries and Theorem~1 to 
rewrite the Hamiltonian as a function on a punctured 
2-dimensional sphere, satisfying
the condition of divergence near the punctures.  Applying Lemma
\ref{lemma:punctures} and recalling that a phase point in a
Hamiltonian system moves along a level set of the Hamiltonian, we
shall be home.

\begin{proof}
  The center of vorticity of the group $\Gamma, \Gamma'$ is at $c$,
  that of the singleton group $-\Gamma-\Gamma'$ at $-c$.  By
  Theorem~1, the vector connecting these centers is a local first
  integral.  Hence passing to the quotient by translations, these
  centers may be assumed immobile.  Within the group $\Gamma,
  \Gamma'$, the position of one vortex determines the position of the
  other (it is at a definite ratio of distances across their center).
  Hence the trajectory of the vortex with vorticity $\Gamma$
  determines the trajectories of all 3 vortices up to translation, and
  the hamiltonian $H$ may be regarded as a function of 
$\zeta = \xi + i\eta$ alone {\it as long as the trajectory of
$\zeta$ lies on a single chart} .
  If the vortices $\Gamma, \Gamma'$ are very
  close, they orbit like a binary star around their immobile center
  $c$ within the chart, so that sooner or later ${\rm arg\,}\zeta$
  increases by $2\pi$.  Since $H(\zeta)\to +\infty$ as $\zeta \to 0$,
  for large enough $E\in \R$ the connected component of $\{ \zeta\in
  \C \diagdown 0 \mid |H(\zeta)| > E \}$ surrounding the singularity
  $\zeta = 0$ is topologically a punctured open disk, free of critical
  points of $H$.  (The infimum of such $E$ is the largest of the
  saddle values of $H$.)  The level sets of $H$ on this neighborhood
  are topologically circles, and so every $\zeta$ starting from
  $\zeta(0)\neq 0$ in this neighborhood returns to $\zeta(0)$,
  guaranteeing relative periodicity.
  
  We must deal with the scenario where the trajectory of $\zeta$ 
does not lie on a single chart.  Since
  $\Gamma/\Gamma'\in \Q$, the lowest common multiple $L$ of $2, 1 +
  \Gamma/\Gamma', 1 + \Gamma'/\Gamma$ makes sense.  To define $\zeta$
  on the whole cylinder, we must swell the cylinder to $\C/L\pi\Z$.
  The swollen cylinder $\C/L\pi\Z$ covers the original cylinder
  $\C/2\pi\Z$ and $H$ as a function of $\zeta$ lifts to a
  function on $\C/L\pi\Z \diagdown \{ {\rm singularities} \}$.
  The singularities represent the collisions between
  $$
  \Gamma\sim\Gamma' {\rm ~(front~and~back)},\quad \Gamma\sim
  -\Gamma-\Gamma',\quad \Gamma'\sim - \Gamma-\Gamma'
  $$
  where $|H| \to \infty$~; off the singularities, by
  (\ref{complexh}),
\begin{equation}\label{threefrogs}
e^{2\pi H/\Gamma\Gamma'} =
\dfrac{\left| \sin\left( c + \dfrac{\zeta}{1 + \Gamma/\Gamma'}
                        \right) \right|^{1 + \Gamma/\Gamma'}
\left| \sin\left( c - \dfrac{\zeta}{1 + \Gamma'/\Gamma}
                        \right)\right|^{1 + \Gamma'/\Gamma}}
{|\sin\zeta|}.
\end{equation}
Toward the `ends' $\eta \to \pm\infty$, $|H| \to \infty$ as well.
Topologically $\C/L\pi\Z \diagdown \{ {\rm singularities} \}$ is a
sphere with at least 4 punctures.  (\ref{threefrogs}) shows that $H$
is Morse and $|H| \to \infty$ near each puncture.  By Lemma
\ref{lemma:punctures}, the generic level sets of $H$ are loops,
representing (putting horizontal translation back in) relative
periodic orbits, and there exist values of $\zeta$ representing
relative equilibria as well as separatrices (relative heteroclinic
orbits) connecting relative equilibria.
\end{proof}

\begin{Remark}
  In Theorem~3, relative periodicity when $\Gamma/\Gamma'\notin \Q$ is
  spoilt only for $\zeta(0)$ too large.  For such $\zeta(0)$, the
  orbit is relative quasi-periodic.  Of course, even when
  $\Gamma/\Gamma' \notin \Q$ there are questions that can be settled
  within a chart.  Thus, for 3 vortices with arbitrary vorticities
  that sum to zero, topological reasons imply the existence
of a configuration that forms a relative equilibrium.  
\end{Remark}

\begin{Theorem}
  Let $b\in \R \diagdown \{ 0 \}$.  On a cylinder, consider the
 configuration of $4$ vortices with vorticities $\Gamma, \Gamma',
  -\Gamma', -\Gamma$ ($\Gamma$ and $\Gamma'$ being of the same sign) at
  $$
  ib + \frac{2\Gamma'}{\Gamma + \Gamma'}\zeta,\quad ib -
  \frac{2\Gamma}{\Gamma + \Gamma'}\zeta,\quad -ib -
  \frac{2\Gamma}{\Gamma + \Gamma'}\overline{\zeta},\quad -ib +
  \frac{2\Gamma'}{\Gamma + \Gamma'}\overline{\zeta} \> ;
  $$
  Let $\Gamma/\Gamma'\neq 1$.  Then for a generic choice of the
  initial condition $\zeta(0)$ these vortices form a relative
  periodic orbit, and for isolated choices of $\zeta(0)$ they form
  a relative equilibrium or a separatrix connecting relative equilibria.
  If\/ $\Gamma/\Gamma' = 1$, the same conclusion holds for $\zeta(0)$
  such that $|{\rm \, Im\,}\zeta(0)| < b$ or 
$\pi H(\zeta(0))/\Gamma^2 < \log {\rm sinh\,}b$.
\end{Theorem}

Combined with Proposition of section 2, Theorem 4 gives relative
equilibria and relative periodic orbits of $N = 4n$ vortices for all
$n \geqslant 1$.

\begin{proof}
  As in the proof of Theorem~3, the positions of all 4 vortices are
  determined by those of the ones with vorticities $\Gamma$ and
  $-\Gamma$.  Thanks to a supplementary reflexive symmetry $z\mapsto
  \overline{z}$, the position of $\Gamma$ determines that of
  $-\Gamma$.  This time, after passing to the quotient by
  translations, $H$ is a genuine function on the cylinder $\C/\pi\Z$
  of $\zeta = \xi+i\eta$, $-\pi/2 < \xi \leqslant \pi/2$, with the
  singularities removed.  Off the singularities, by (\ref{complexh}),
\begin{equation}\label{fourfrogs}
e^{2\pi H/\Gamma\Gamma'} =
\left|\dfrac{\sin\left( 
 ib + \dfrac{\Gamma'\zeta+\Gamma\overline{\zeta}}{\Gamma+\Gamma'} 
\right)}
{\sin\zeta} \right|^2
\left| \sin\left(
 ib + \dfrac{\zeta-\overline{\zeta}}{1 + \Gamma/\Gamma'}
\right) \right|^{\Gamma/\Gamma'}
\left| \sin\left(
 ib - \dfrac{\zeta-\overline{\zeta}}{1 + \Gamma'/\Gamma}
\right) \right|^{\Gamma'/\Gamma}.
\end{equation}
\noindent In particular, when $\Gamma/\Gamma' = 1$,
\begin{equation}\label{realfourfrogs}
e^{2\pi H/\Gamma^2} =
\frac{\sin^2\xi + {\rm sinh}^2b}{\sin^2\xi + {\rm sinh}^2\eta}
\left| {\rm \, sinh\,}(b + \eta){\rm \, sinh\,}(b - \eta) \right|.
\end{equation}
\noindent
The isolated singularities represent simultaneous collisions between
$$
\Gamma\sim\Gamma' {\rm ~and~} -\Gamma'\sim -\Gamma
$$
where $H \to +\infty$, and, if $\Gamma/\Gamma'\neq 1$, between
$$
\Gamma\sim -\Gamma' {\rm ~and~} \Gamma'\sim -\Gamma
$$
where $H \to -\infty$.  Toward the ends, $H\to +\infty$.  There are
also circles of singularities $\eta = -b(1 + \Gamma/\Gamma')/2$, $b(1
+ \Gamma'/\Gamma)/2$ representing collisions between
$$
\Gamma\sim -\Gamma,\quad \Gamma'\sim -\Gamma'
$$
where $H \to -\infty$.  Let us saw the cylinder $\C/\pi\Z$ of
$\zeta$ into 3 trunks:
\begin{equation*}
\begin{split}
  &C_+ = \{ \zeta \> | \> b(1 + \Gamma'/\Gamma)/2  < \eta \}, \\
  &C_0 = \{ \zeta \> | \>
  -b(1 + \Gamma/\Gamma')/2 < \eta < b(1 + \Gamma'/\Gamma)/2 \}, \\
  &C_- = \{ \zeta \> | \> \eta < -b(1 + \Gamma/\Gamma')/2 \}.
\end{split}
\end{equation*}
Topologically $C_+, C_0, C_-$ are spheres with punctures.  $C_0$
contains $\zeta = 0$, the simultaneous collisions between
$\Gamma\sim\Gamma', -\Gamma'\sim -\Gamma$, so $C_0$ has at least 3
punctures and $|H|\to \infty$ near each puncture.  Lemma
\ref{lemma:punctures} applies to $C_0$ and implies the existence of
relative periodic orbits and relative equilibria.

For the moment, suppose $\Gamma/\Gamma'\neq 1$.  $\zeta$ representing
the simultaneous collisions between $\Gamma\sim -\Gamma', \Gamma'\sim
-\Gamma$ is in $C_+$ or $C_-$ accordingly as $\Gamma/\Gamma' > 1$ or
$< 1$.  If $\Gamma/\Gamma' > 1$, this puts on $C_+$ at least 3
punctures near each of which $|H| \to \infty$, so Lemma
\ref{lemma:punctures} applies and implies the existence of relative
periodic orbits and relative equilibria; whereas $C_-$ acquires only 2
punctures, so we can conclude the existence only of relative periodic
orbits.  If $\Gamma/\Gamma' < 1$, the r\^oles of $C_+, C_-$ are
reversed.

Note that as $H$ is symmetric under the lateral reflection along $\xi
= 0$ and along $\xi = \pi/2$, every point on either line where
$\partial H/\partial \eta$ vanishes is critical.  Let $\Gamma/\Gamma'
> 1$ and work on $C_+$.  The strip $0 < \xi < \pi/2$ is free of
critical points, for here by (\ref{realfourfrogs}) $H$ is strictly
monotone in $\xi$ along any line $\eta =$ constant.  Along $\xi = 0$,
$H\to -\infty$ as $\eta \to b(1 + \Gamma'/\Gamma)/2$ or $b(\Gamma +
\Gamma')/(\Gamma - \Gamma')$, between which $\partial H/\partial \eta$
must vanish, signaling a saddle at say $\zeta_1$.  Along $\xi =
\pi/2$, $H\to -\infty$ or $+\infty$ as $\eta \to b(1 +
\Gamma'/\Gamma)/2$ or $+\infty$.  These bits of information, together
with the fact that $H$ is Morse, imply that $\partial H/\partial \eta$
vanishes twice along $\xi = \pi/2$, signaling a maximum at say
$\zeta_2$ and a saddle (which shall be left nameless).  As a bonus we
learn that 2 relative equilibria are represented in $C_+$, whereas a
count of 3 singularities just predicts at least 1 relative
equilibrium.  The analysis works mutatis mutandis on $C_-$ if
$\Gamma/\Gamma' < 1$.

Finally, suppose $\Gamma/\Gamma' = 1$.  Then the simultaneous
collisions $\Gamma\sim -\Gamma', \Gamma'\sim -\Gamma$ as well as
$\zeta_1, \zeta_2$ escape to the ends $\eta\to \pm\infty$, and toward
the ends $2\pi H/\Gamma^2$ asymptotes to $\log (\sin^2\xi + {\rm
  sinh}^2b)$, which remains bounded.  Hence all the critical points in
$C_+, C_-$ disappear.  Relative periodic orbits are represented by
compact level sets of $H$, i.e.\ those that fill the region $e^{\pi
  H/\Gamma^2} < {\rm sinh\,}b$ of $C_+, C_-$~; there is no relative
equilibrium on these trunks.
\end{proof}

\begin{figure*}[ht]
  \begin{center}
    \subfigure {\label{} \includegraphics[angle=0,
      width=5cm]{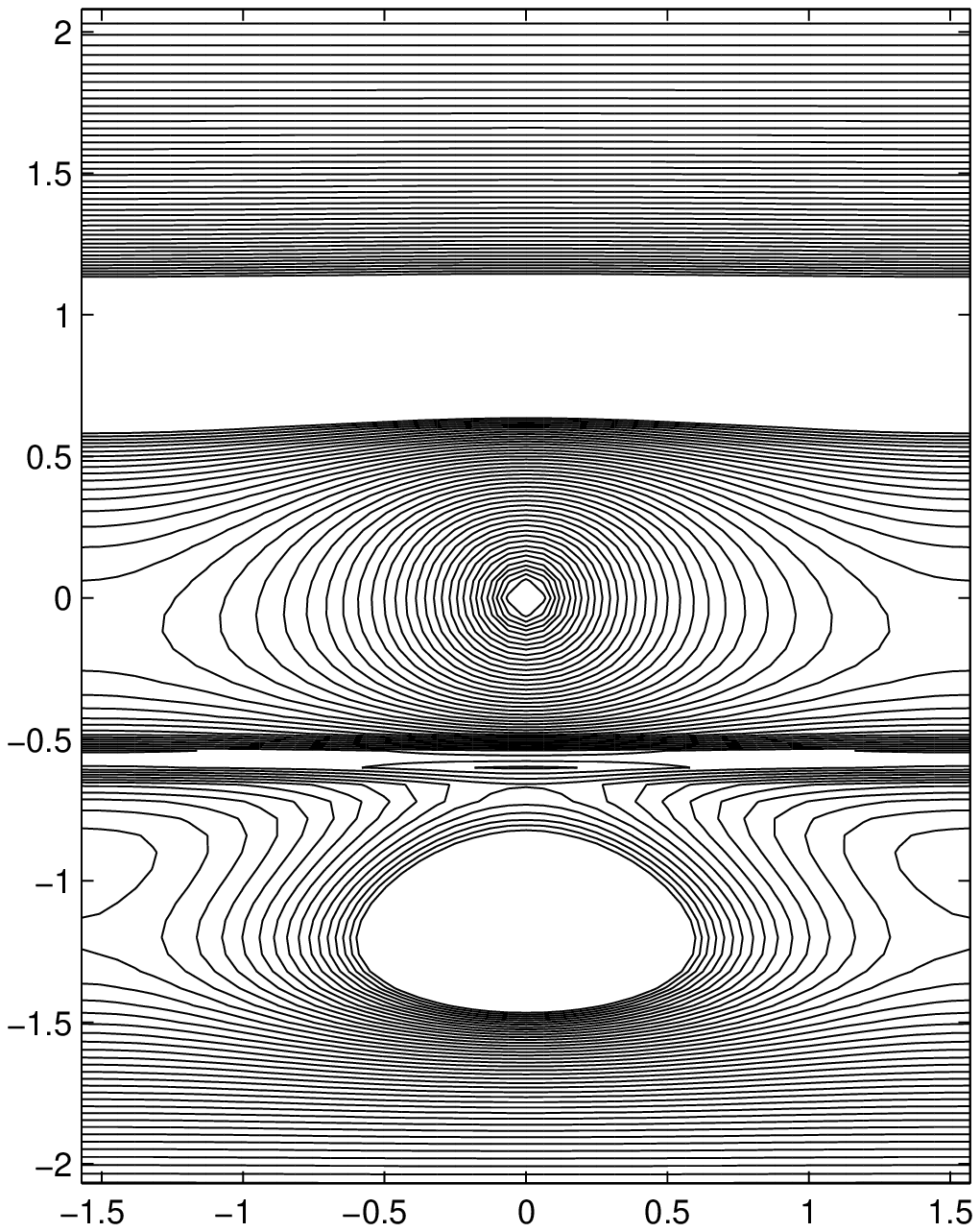}} \subfigure {\label{}
      \includegraphics[angle=0, width=5cm]{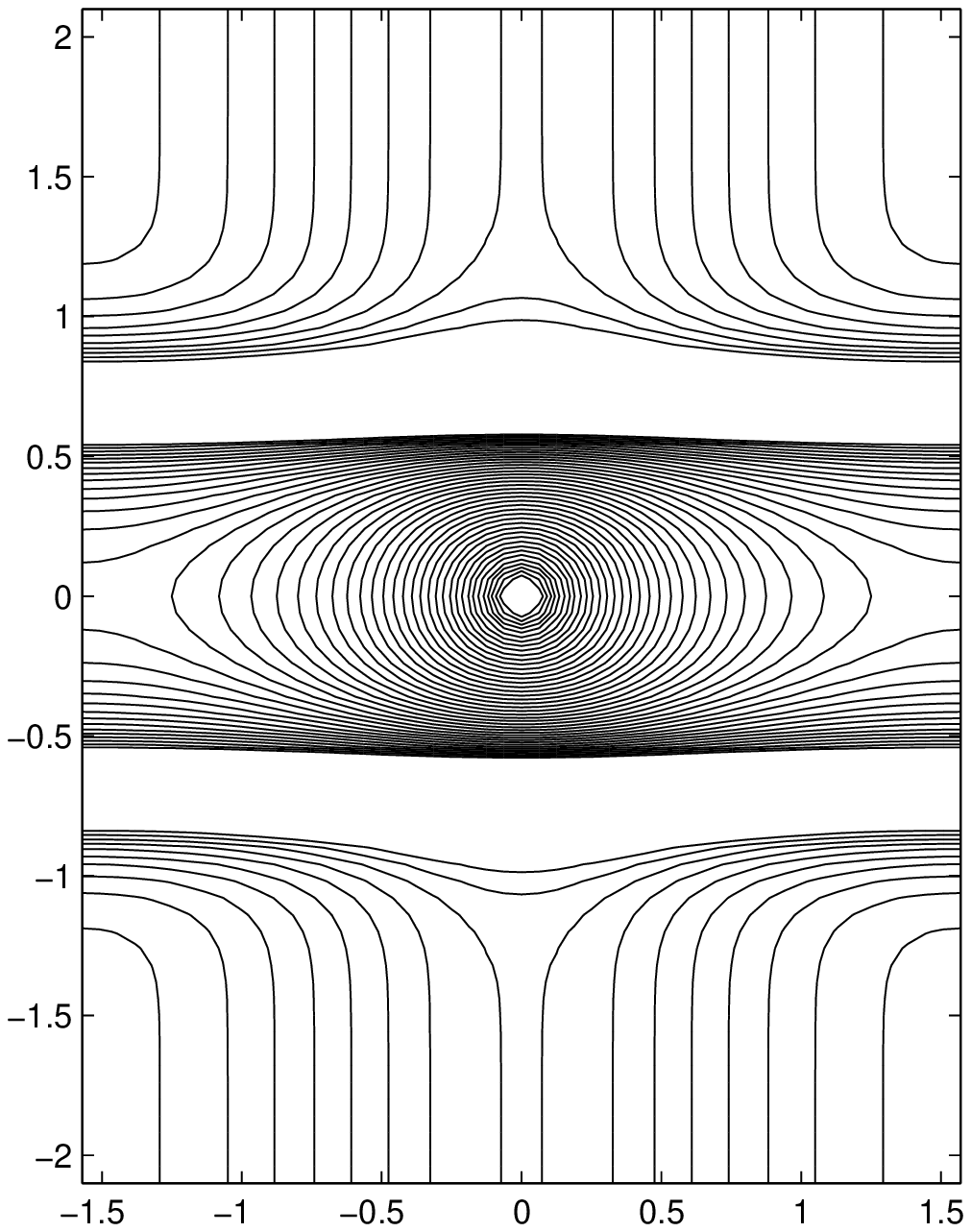}} \subfigure
    {\label{} \includegraphics[angle=0, width=5cm]{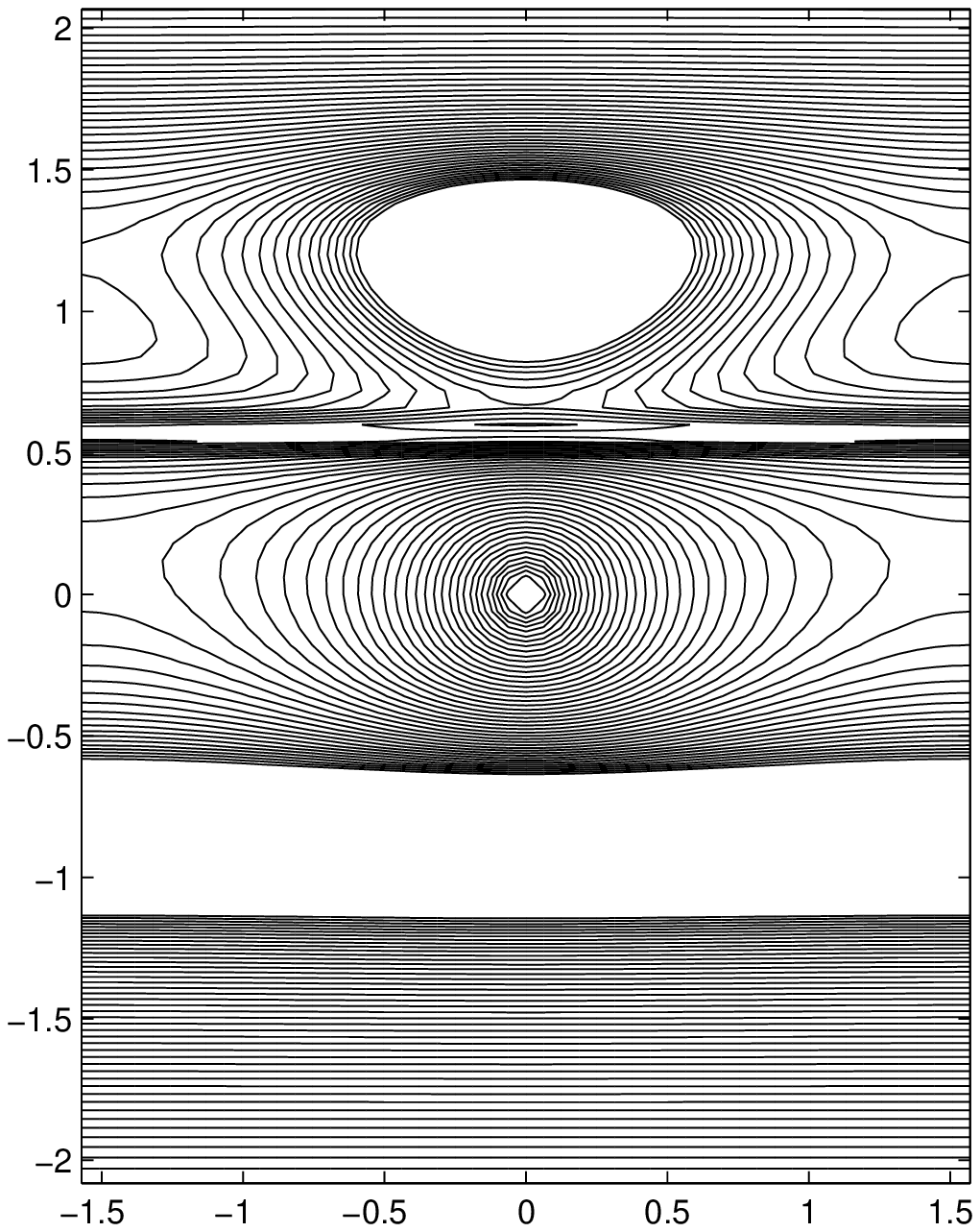}}
    \caption{}
    \end{center}
\end{figure*}

The plots of fig.~4 depict the level sets of $H$ as a function
of $\zeta$ for $\Gamma/\Gamma' <
1$, $= 1$, $> 1$ respectively; they were drawn at $b = 1$.  By
(\ref{fourfrogs}), the levels for $\Gamma/\Gamma' < 1$ and $> 1$ are
mirror images of each other via $\zeta \mapsto \overline{\zeta}$.  The
blank holes and bands indicate where $H$ diverges to $-\infty$ too
steeply, while the diamond in the middle of each plot surrounds a peak
$H \to +\infty$.

Take the $N = 4$ case as in Theorem 4 and initially align the 4
vortices vertically:
$\xi(0) = 0$.  If $\eta(0)$ is sufficiently small, the vortices of the
group $\Gamma, \Gamma'$ orbit like a binary counter-clockwise, the
vortices of the group $-\Gamma', -\Gamma$ orbit like a binary
clockwise, while the 2 groups progress together like a vortex pair.
The superposition produces {\it leapfrogging\/}, a relative periodic
orbit whose plane limit $r\to \infty$ is observed as the motion of a
cross-section of consecutive vortex rings as they overtake each other.
By adjusting the parameters $\Gamma/\Gamma'$, $b$, $\zeta(0)$, we can
render leapfrogging on a cylinder not only relative periodic but
periodic.  Alternatively, if $\eta(0)$ is sufficiently close to $b(1 +
\Gamma'/\Gamma)/2$ or to $-b(1 + \Gamma/\Gamma')/2$, the vortices
$\Gamma', -\Gamma'$ or $\Gamma, -\Gamma$ form a pair and rush off
without leapfrogging.  In the planar theory, in the case
$\Gamma/\Gamma' = 1$, \cite{Love} calculated the critical value of
$\eta(0)$ that separates the leapfrogging and non-leapfrogging
r\'egimes.  In our setup this value may be obtained at
once as follows.  

In the situation of Theorem~4, denote by $\rho(b, \Gamma/\Gamma')$
the distance from the origin $\zeta = 0$ to the nearest
separatrix.  Then $\eta(0) = \rho(b, \Gamma/\Gamma) = \rho(b, 1)$.
Denote by $\zeta_{\rm re} = \xi_{\rm re} + i\eta_{\rm re}$ a value of
$\zeta$ at a saddle of $H(\zeta)$, representing a relative
equilibrium.  Inside the separatrices connecting the saddles we have
leapfrogging; outside, not.  $\rho =\rho(b, 1)$ is the ordinate at which
a separatrix cuts the $\eta$-axis.  Since the value of $H$ is the same
along the separatrices as on the saddles, $H(0, \rho) = H(\xi_{\rm
  re}, \eta_{\rm re})$.  It is clear that a relative equilibrium
occurs when 2 vortex pairs are antipodal: $\xi_{\rm re} = \pm\pi/2,
\eta_{\rm re} = 0$.  This fixes $\rho$ in the cylindrical theory:
$\sqrt{2}{\rm \, tanh\,}\rho = {\rm tanh\,}b$.  Restoring $r$ and taking
the plane limit $r\to \infty$, we get in the planar theory $\rho =
b/\sqrt{2}$, agreeing with \cite[section~3]{Love}, which arrived at
$(b + \rho)/(b - \rho) = 3 + 2\sqrt{2}$.

When $\Gamma/\Gamma'\neq 1$, $\zeta_{\rm re}$ and $\rho(b,
\Gamma/\Gamma')$ are difficult to pin down in closed form.  At any
rate $\xi_{\rm re} = \pm\pi/2$~; $\eta_{\rm re}$ is the unique root of
\begin{eqnarray*}
(\Gamma + \Gamma'){\rm \, tanh\,}\eta \! &+&\! (\Gamma - \Gamma'){\rm \,
  tanh}\left( b - \frac{\Gamma - \Gamma'}{\Gamma + \Gamma'}\eta
\right)\\
\! &-&\! \Gamma{\rm \, coth}\left( b +\frac{2\eta}{1 + \Gamma/\Gamma'}\right )
+\Gamma'{\rm \, coth}\left( b -\frac{2\eta}{1 + \Gamma'/\Gamma}\right)
= 0
\end{eqnarray*}
which in view of (\ref{realeq}) is the condition that the
vertically aligned pairs $\Gamma, -\Gamma$ and $\Gamma', -\Gamma'$,
antipodal to each other, move with the same velocity.  If
$\Gamma/\Gamma' = 1 + \epsilon$, then up to 2nd order in $\epsilon$,
$$
\eta_{\rm re} \simeq {\rm tanh\,}b\, {\rm \, sech}^2b \left(
  \frac{\epsilon}{2} - \left( 1 + \frac{{\rm sech}^4b}{2}
  \right)\frac{\epsilon^2}{4} \right), \qquad \rho(b, 1 + \epsilon) =
\rho(b, 1) - \frac{{\rm tanh\,}b\, {\rm \, sech}^2b}{1 + {\rm cosh}^2b} \,
\frac{\epsilon^2}{4\sqrt{2}}.
$$

\begin{Remark}
  By an argument parallel to that of Theorem~4 we see that $4$
  vortices with vorticities $\Gamma, \Gamma', -\Gamma', -\Gamma$ at
  $$
  a + \frac{2\Gamma'}{\Gamma + \Gamma'}\zeta,\quad a -
  \frac{2\Gamma}{\Gamma + \Gamma'}\zeta,\quad -a +
  \frac{2\Gamma}{\Gamma + \Gamma'}\overline{\zeta},\quad -a -
  \frac{2\Gamma'}{\Gamma + \Gamma'}\overline{\zeta}
  $$
  leapfrog as well (fig.~5, left diagram).
\begin{figure}[ht]
  \begin{center}
    \includegraphics[width=12cm]{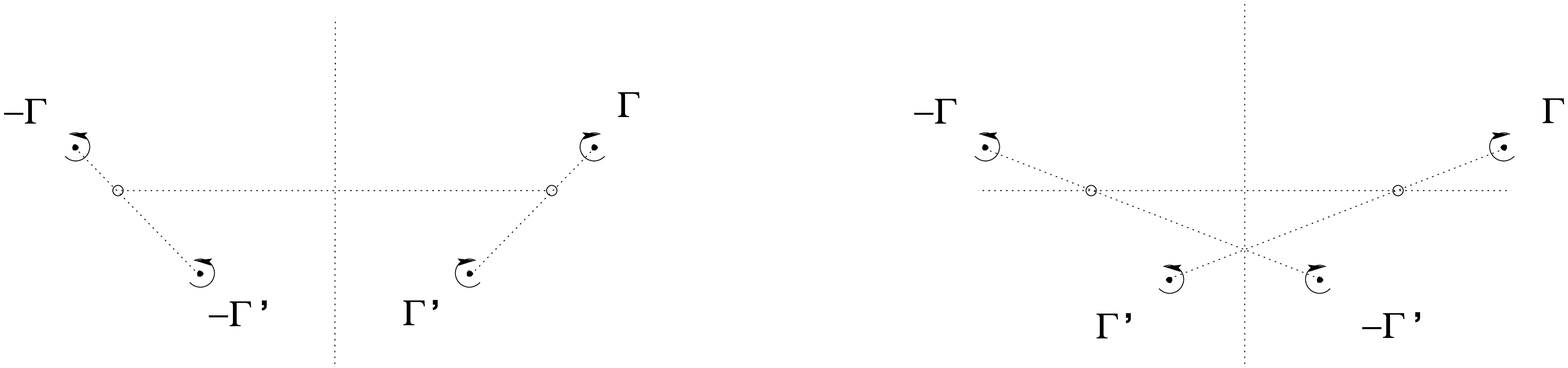}
    \caption{}
  \end{center}
\end{figure}

\noindent 
Unlike the $N = 4$ case of Theorem~4, however, the configuration on
the right does not leapfrog.
\end{Remark}

\begin{Remark}
  Leapfrogging vortices and their generalizations analysed above owe
  their relative periodicity to the type of symmetry compatible with
  the local first integral of Theorem~1.  Other types of symmetry
  permit other types of relative periodic orbits.  Thus, $2n$ vortices
  with identical vorticity $\Gamma$ at (\ref{uppertrain}),
  (\ref{lowertrain}) form a relative periodic orbit
  \cite[section~3.2]{SouliereTokieda}.
\end{Remark}

\begin{Remark}
  Vortex streets and leapfrogging vortices can be adapted to a torus,
  where they form relative periodic orbits.  A torus accommodates many
  further types of relative periodic orbits.  For example on
  $\C/(\pi\Z + i\pi\Z)$, by splitting each point of a sub-lattice
  into a rectangular quadruplet of vortices with vorticities $\Gamma,
  -\Gamma, \Gamma, -\Gamma$, we create a periodic orbit, the `dancing
  vortices' of \cite{Tokieda}.
\end{Remark}

\noindent {\bf Acknowledgments}.
The work of JM was partially supported by the
European Union through the Research Training Network MASIE.  
TT thanks Yiannis Petridis and Morikazu Toda for instructive
conversations and David Acheson for his gift of a copy of
the paper \cite{Love}.

\end{document}